\documentclass[11pt]{article}
\usepackage{amsmath}
\usepackage{amssymb, amscd, amsmath, amsthm, amsfonts, ae}
\usepackage{graphicx,array}
\newcommand{\bc}{\mathbb{C}}
\newcommand{\bz}{\mathbb{Z}}

\newcommand{\p}{\frac{\partial}}
\newcommand{\x}{\partial x}

\title{A Class of Irreducible Modules for the Extended Affine Lie Algebra
$\widetilde{\frak{gl}_l({\mathbb{C}_q})}$}
\author{Ziting Zeng\footnote{The author gratefully
acknowledges the grant support from NNSF of China No.10801010.}
\\School of Mathematical Science, Beijing Normal University \\
Beijing, China 100875\\zengzt@bnu.edu.cn}
\date{}
\begin{document}
\maketitle
\begin{abstract}
We construct a class of modules for extended affine Lie algebra
$\widetilde{\frak{gl}_l({\bc_q})}$ by using the free fields. A necessary and sufficient
condition is given for those modules being irreducible.

\end{abstract}

\newtheorem{theorem}{Theorem}[section]
\newtheorem{remark}{Remark}[section]
\newtheorem{corollary}{Corollary}[section]
\allowdisplaybreaks

\section{Introduction}
The representation theory of affine Kac-Moody Lie algebras has
remarkable applications in many areas of mathematics and
mathematical physics.
Extended affine Lie algebras are a higher dimensional generalization
of affine Kac-Moody Lie algebras first introduced by [H-KT] under the name of quasi-simple Lie algebras and
systematically studied in [AABGP]. This family of  newly developed Lie algebras includes
toroidal Lie algebras and the central extensions of the matrix Lie algebras coordinated by quantum tori as examples.
The study of representations
for extended affine Lie algebras coordinated by quantum tori have drawn a lot of attentions recently via various module realizations. For instance,
the representations via vertex operator construction was obtained in [BS], [G1], [G2], [BGT]. Fermionic and bosonic realizations
were constructed in [G3] and [L]. Unitary representations were studied in [JK], [ER], [GZ1] and [Z]. A deformed version was given in [VV].

The Wakimoto's  free fields construction provides a remarkable way
to realize  affine Kac-Moody Lie algebras (see [W1], [FF] and [W2]). In [GZ] and [Z], this approach  has been successfully used to construct the
representation of $\widetilde{\frak{gl}_l({\mathbb{C}_q})}$, the
extended affine Lie algebra of type $A$ coordinated by a quantum
torus, for $l=2,3$. Motivated by [Z] and [GZ2], we shall
construct in this paper a class of irreducible representations of the extended
affine Lie algebra $\widetilde{\frak{gl}_l(\bc_q)}$ which contains
$l=2$ in [GZ] as its special case. We further give a necessary and sufficient
condition for those modules being irreducible.

 The paper is
organized as follows.  In Sect 2, we recall some
notations and definitions. Our main theorem will appears
in Sect 3 where we construct a class of modules for
$\widetilde{\frak{gl}_l(\bc_q)}$ by using the free fields.  Then we work out
 the condition of
irreducibility for modules constructed in Sect 4.

\section{Preliminary}

 Let $q$ be a non-zero complex number. A
quantum $2$-torus associated to $q$ is the unital associative
$\bc$-algebra $\bc_{q}[s^{\pm 1}, t^{\pm 1}]$ (or, simply $\bc_{q}$)
with generators $s^{\pm 1},  t^{\pm 1}$ and relations
$$s s^{-1} =s^{-1}s =t t^{-1}= t^{-1}t=1 \, \text{  and } \, \ ts = q st.$$

Let $d_s$, $d_t$ be the degree operators on $\bc_q$ defined by
$$d_s(s^mt^n) = ms^m t^n, \, \,  d_t(s^m t^n) = n s^m t^n $$
 for $m, n\in \bz$.

Let $\frak{gl}_l$ be the general linear algebra, which is one
dimensional central extension of the simple Lie algebra $\frak{sl}_l$ of type
$A_{l-1}$. Denote $\dot{\frak{h}}$ the Cartan
subalgebra of $\frak{gl}_l$ which is the set of the diagonal
matrices, and $\dot{P}\subset \frak{\dot{\frak{h}}}^*$ be the set of
weight lattice. Its Borel subalgebra denoted as $\dot{\frak{b}}$ is
the set of all upper triangular matrices.

We form a  natural central extension of $\frak{gl}_{l}\otimes\bc_q$
as follows.
$$\widehat{\frak{gl}_{l}(\bc_q)} = \frak{gl}_{l}\otimes\bc_q\oplus\bc c_s\oplus
\bc c_t $$ with Lie bracket
\begin{align*} & [E_{ij}\otimes s^{m_1}t^{n_1}, E_{kn}\otimes s^{m_2}t^{n_2}]\\
=&\delta_{jk}q^{n_1m_2}E_{in}\otimes
s^{m_1+m_2}t^{n_1+n_2}-\delta_{in}q^{n_2 m_1}E_{kj}
\otimes s^{m_1+m_2}t^{n_1+n_2}\\
&+m_1q^{n_1 m_2}\delta_{jk}\delta_{in}\delta_{m_1+m_2, 0}\delta_{n_1
+ n_2, 0}c_s +n_1q^{n_1m_2}\delta_{jk}\delta_{in}\delta_{m_1+m_2,
0}\delta_{n_1+n_2, 0}c_t
\end{align*}
for $m_1, m_2, n_1, n_2\in \bz$, $1\leq i, j, k, n\leq l$, where
$E_{ij}$ is the matrix whose $(i, j)$-entry is $1$ and $0$
elsewhere, and $c_s$ and $c_t$ are central elements of
$\widehat{\mathfrak{gl}_{l}(\bc_q)}$.

The derivations $d_s$ and $d_t$ can be extended to derivations on
$\widehat{\frak{gl}_{l}(\bc_q)}$. Now we can define the semi-direct
product of the Lie algebra $\widehat{\frak{gl}_{l}(\bc_q)}$ and
those derivations:
$$\widetilde{\frak{gl}_{l}(\bc_q)} =\widehat{\frak{gl}_{l}(\bc_q)} \oplus \bc d_s
\oplus \bc d_t. $$

The Lie algebra $\widetilde{\frak{gl}_{l}(\bc_q)}$ is  an extended
affine Lie algebra of type $A_{l-1}$ with nullity $2$. ( See [AABGP] and [BGK]
for definitions).

So the Cartan subalgebra $\frak{h}$ of
$\widetilde{\frak{gl}_{l}(\bc_q)}$ is
$$\frak{h}=\dot{\frak{h}}\oplus \bc c_s\oplus \bc c_t\oplus \bc d_s\oplus \bc d_t,$$
the Borel subalgebras $\frak{b}$ is
$$\frak{b}=\dot{\frak{b}}\otimes \bc_q \oplus \bc c_s\oplus \bc c_t\oplus \bc d_s\oplus \bc d_t$$
 and denote
$\frak{n}_+=\bc\left\{E_{ij}\otimes \bc_q,j>i\geq 1\right\}$ the
nilpotent radical of $\frak{b}$.

\section{The module}
In this section, we will construct the modules of
$\widetilde{\frak{gl}_{l}(\bc_q)}$ by using the free fields.

 Let
 $$V=\bc [x_i(m,n),
i=2,\cdots,l, m,n\in\bz]$$ be a polynomial ring with infinitely many
variables $x_i(m,n)$, we define the following operators:

\begin{eqnarray*}
e_{i1}(m_1,n_1)&=&x_i(m_1,n_1),  \  i=2,\cdots, l ;\\
e_{1i}(m_1,n_1)&=&q^{-m_1n_1}\mu\p{\x_i(-m_1,-n_1)}\\&
&-\sum_{(m,n)\in \bz^2\atop (m',n')\in \bz^2}
q^{n_1m'+nm_1+nm'}x_2(m_1+m+m',n_1+n+n')\p{\x_2(m,n)}\p{\x_i(m',n')}\\
& &-\sum_{(m,n)\in \bz^2\atop (m',n')\in \bz^2}
q^{n_1m'+nm_1+nm'}x_3(m_1+m+m',n_1+n+n')\p{\x_3(m,n)}\p{\x_i(m',n')}\\
& &-\cdots\\
& &-\sum_{(m,n)\in \bz^2\atop (m',n')\in \bz^2}
q^{n_1m'+nm_1+nm'}x_l(m_1+m+m',n_1+n+n')\p{\x_l(m,n)}\p{\x_i(m',n')},\\
&  &\text{ for } i=2,\cdots, l;\\
e_{ij}(m_1,n_1)&=&\sum_{(m,n)\in \bz^2} q^{mn_1}x_i(m_1+m,n_1+n)\p{\x_j(m,n)}, \  i,j=2,\cdots, l;\\
e_{11}(m_1,n_1)&=&\mu\delta_{(m_1,n_1),(0,0)}-\sum_{(m,n)\in \bz^2} q^{nm_1}x_2(m_1+m,n_1+n)\p{\x_2(m,n)}\\
& &-\sum_{(m,n)\in \bz^2} q^{nm_1}x_3(m_1+m,n_1+n)\p{\x_3(m,n)}\\
& &-\cdots\\
& &-\sum_{(m,n)\in \bz^2} q^{nm_1}x_l(m_1+m,n_1+n)\p{\x_l(m,n)}\\
 D_1&=&\sum_{i=2}^{l}\sum_{(m,n)\in \bz^2}mx_i(m,n)\p{\x_i(m,n)}\\
 D_2&=&\sum_{i=2}^{l}\sum_{(m,n)\in \bz^2}nx_i(m,n)\p{\x_i(m,n)}
\end{eqnarray*}

Though some operators are infinite sums, they are well-defined when acting on $V$ (as only
 finite summands left). Here we give our first theorem:

\begin{theorem} There is a Lie algebra homomorphism $\phi:
\widetilde{\frak{gl}_l({\mathbb{C}_q})} \rightarrow \frak{gl}(V)$ given by
$$\phi (E_{ij}\otimes s^m t^n)=e_{ij}(m,n),i,j=1,\cdots l, $$
$$\phi (d_s)=D_1,\phi(d_t)=D_2,\phi(c_s)=\phi(c_t)=0.$$
Hence $V$ is a module of $\widetilde{\frak{gl}_l(\bc_q)}$.
\end{theorem}
\noindent{\sl Proof}:\quad We only need to check that $\phi$ preserves the brackets.
 The following formidable
calculation is to check them case by case. For neatness we simplify
$\sum_{(m,n)\in \bz^2\atop (m',n')\in \bz^2}$ or $\sum_{(m,n)\in
\bz^2}$ as $\sum$.

 If $i\neq j$
\begin{align*} &[e_{i1}(m_1,n_1),e_{1j}(m_2,n_2)]\\
=&[x_i(m_1,n_1),-\sum
q^{n_2m'+nm_2+nm'}x_i(m_2+m+m',n_2+n+n')\p{\x_i(m,n)}\p{\x_j(m',n')}]\\
=&\sum
q^{n_2m'+n_1m_2+n_1m'}x_i(m_1+m_2+m',n_1+n_2+n')\p{\x_j(m',n')}\\
=&q^{n_1m_2}e_{ij}(m_1+m_2,n_1+n_2).\end{align*}

\begin{align*} &[e_{i1}(m_1,n_1),e_{1i}(m_2,n_2)]\\
=&[x_i(m_1,n_1),q^{-m_2n_2}\mu\p{\x_i(-m_2,-n_2)}]\\
&+[x_i(m_1,n_1),-\sum
q^{n_2m'+nm_2+nm'}x_2(m_2+m+m',n_2+n+n')\p{\x_2(m,n)}\p{\x_i(m',n')}]\\
 &+\cdots\\
&+[x_i(m_1,n_1),-\sum
q^{n_2m'+nm_2+nm'}x_i(m_2+m+m',n_2+n+n')\p{\x_i(m,n)}\p{\x_i(m',n')}]\\
&+\cdots\\
&+[x_i(m_1,n_1),-\sum
q^{n_2m'+nm_2+nm'}x_l(m_2+m+m',n_2+n+n')\p{\x_l(m,n)}\p{\x_i(m',n')}]\\
=&-q^{-m_2n_2}\mu\delta_{(m_1+m_2,n_1+n_2),(0,0)}+\sum q^{n_2m_1+nm_2+nm_1}x_2(m_2+m+m_1,n_2+n+n_1)\p{\x_2(m,n)}\\
&+\cdots\\
&+\sum q^{n_2m_1+nm_2+nm_1}x_i(m_2+m+m_1,n_2+n+n_1)\p{\x_i(m,n)}\\
&+\sum q^{n_2m'+n_1m_2+n_1m'}x_i(m_2+m_1+m',n_2+n_1+n')\p{\x_i(m',n')}\\
&+\cdots\\
&+\sum q^{n_2m_1+nm_2+nm_1}x_l(m_2+m+m_1,n_2+n+n_1)\p{\x_l(m,n)}\\
=&-q^{m_1n_2}\mu\delta_{(m_1+m_2,n_1+n_2),(0,0)}+q^{n_2m_1}(\sum
q^{n(m_2+m_1)}x_2(m_2+m+m_1,n_2+n+n_1)\p{\x_2(m,n)}\\
&\quad\quad \cdots\\
&\quad\quad +\sum
q^{n(m_2+m_1)}x_l(m_2+m+m_1,n_2+n+n_1)\p{\x_l(m,n)}\\
&+q^{n_1m_2}\sum
q^{(n_2+n_1)m'}x_i(m_2+m_1+m',n_2+n_1+n')\p{\x_i(m',n')}\\
=&-q^{n_2m_1}e_{11}(m_1+m_2,n_1+n_2)+q^{n_1m_2}e_{ii}(m_1+m_2,n_1+n_2).
\end{align*}

\begin{eqnarray*} &[e_{i1}(m_1,n_1),e_{jk}(m_2,n_2)]\\
=&[x_i(m_1,n_1),\sum
q^{mn_2}x_j(m_2+m,n_2+n)\p{\x_k(m,n)}]\\
=&-\delta_{ki} q^{m_1n_2}x_j(m_1+m_2,n_1+n_2)\\
=&-\delta_{ki} q^{m_1n_2}e_{j1}(m_1+m_2,n_1+n_2).\end{eqnarray*}

\begin{eqnarray*} &[e_{i1}(m_1,n_1),e_{11}(m_2,n_2)]\\
=&[x_i(m_1,n_1),-\sum
q^{nm_2}x_i(m_2+m,n_2+n)\p{\x_i(m,n)}]\\
=&q^{n_1m_2}x_i(m_1+m_2,n_1+n_2)\\
=&q^{n_1m_2}e_{i1}(m_1+m_2,n_1+n_2).\end{eqnarray*}

\begin{align*} &[e_{1i}(m_1,n_1),e_{1i}(m_2,n_2)]\\
=&[q^{-m_1n_1}\mu\p{\x_i(-m_1,-n_1)}\\&-\sum q^{n_1m'+nm_1+nm'}
x_j(m_1+m+m',n_1+n+n')\p{\x_j(m,n)}\p{x_i(m',n')},\\
&q^{-m_2n_2}\mu\p{\x_i(-m_2,-n_2)}\\&-\sum q^{n_2m'+nm_2+nm'}
x_j(m_2+m+m',n_2+n+n')\p{\x_j(m,n)}\p{x_i(m',n')}]\\
=&
-\mu q^{m_1n_1}\sum_{{-m_1=m_2+m+m'}\atop{-n_1=n_2+n+n'}}q^{n_2m'+nm_2+nm'}\p{\x_i(m,n)}\p{\x_i(m',n')}\\
&+\mu q^{-m_2n_2}\sum_{{-m_2=m_1+m+m'}\atop{-n_2=n_1+n+n'}}q^{n_1m'+nm_1+nm'}\p{\x_i(m,n)}\p{\x_i(m',n')}\\
&+\sum
q^{n_1(m_2+m''+m')+nm_1+n(m_2+m''+m')+n_2m'+n''m_2+n''m'}\\
&\cdot x_j(m_1+m+m_2+m''+m',n_1+n+n_2+n''+n')
\p{\x_j(m,n)}\p{\x_i(m'',n'')}\p{\x_i(m',n')}\\
&+\sum
q^{n_1m'+(n_2+n+n'')m_1+(n_2+n+n'')m'+n_2m''+nm_2+nm''}\\
&\cdot x_j(m_1+m+m_2+m''+m',n_1+n+n_2+n''+n')
\p{\x_j(m,n)}\p{\x_i(m'',n'')}\p{\x_i(m',n')}\\
&-\sum
q^{n_2(m_1+m''+m')+nm_2+n(m_1+m''+m')+n_1m'+n''m_1+n''m'}\\
&\cdot x_j(m_1+m+m_2+m''+m',n_1+n+n_2+n''+n')
\p{\x_j(m,n)}\p{\x_i(m'',n'')}\p{\x_i(m',n')}\\
&-\sum q^{n_2m'+(n_1+n+n'')m_2+(n_1+n+n'')m'+n_1m''+nm_1+nm''}\\
&\cdot x_j(m_1+m+m_2+m''+m',n_1+n+n_2+n''+n')
\p{\x_j(m,n)}\p{\x_i(m'',n'')}\p{\x_i(m',n')}\\
=&0.\end{align*}

If $i\neq j$,
\begin{align*} &[e_{1i}(m_1,n_1),e_{1j}(m_2,n_2)]\\
=&[-\sum
q^{n_1m'+nm_1+nm'}x_k(m_1+m+m',n_1+n+n')\p{\x_k(m,n)}\p{\x_i(m',n')},\\
&-\sum
q^{n_2m'+nm_2+nm'}x_k(m_2+m+m',n_2+n+n')\p{\x_k(m,n)}\p{\x_j(m',n')}]\\
&+[q^{-m_1n_1}\mu\p{\x_i(-m_1,-n_1)},\\&-\sum
q^{n_2m'+nm_2+nm'}x_i(m_2+m+m',n_2+n+n')\p{\x_i(m,n)}\p{\x_j(m',n')}]\\
&+[-\sum q^{n_1m'+nm_1+nm'}x_j(m_1+m+m',n_1+n+n') \p{\x_j(m,n)}
\p{\x_i(m',n')},\\&q^{-m_2n_2}\mu\p{\x_j(-m_2,-n_2)}]\\
=&\sum
q^{n_1(m_2+m''+m')+nm_1+n(m_2+m''+m')+n_2m'+n''m_2+n''m'}\\
&\cdot
x_k(m_1+m_2+m+m'+m'',n_1+n_2+n+n'+n'')\p{\x_k(m,n)}\p{\x_i(m'',n'')}\p{\x_j(m',n')}\\
&+\sum q^{n_1m'+(n_2+n+n'')m_1+(n_2+n+n'')m'+n_2m''+nm_2+nm''}\\
&\cdot x_k(m_1+m_2+m+m'+m'',n_1+n_2+n+n'+n'')\p{\x_i(m',n')}\p{\x_k(m,n)}\p{\x_j(m'',n'')}\\
&-\sum q^{n_2(m_1+m'+m'')+nm_2+n(m_1+m'+m'')+n_1m'+n''m_1+n''m'}\\
&\cdot x_k(m_1+m_2+m+m'+m'',n_1+n_2+n+n'+n'')\p{\x_k(m,n)}\p{\x_j(m'',n'')}\p{\x_i(m',n')}\\
&-\sum q^{n_2m'+(n_1+n+n'')m_2+(n_1+n+n'')m'+n_1m''+nm_1+nm''}\\
&\cdot x_k(m_1+m_2+m+m'+m'',n_1+n_2+n+n'+n'')\p{\x_j(m',n')}\p{\x_k(m,n)}\p{\x_i(m'',n'')}\\
&-\mu q^{-m_1n_1}\sum_{{-m_1=m_2+m+m'}\atop{-n_1=n_2+n+n'}}q^{n_2m'+nm_2+nm'}\p{\x_i(m,n)}\p{\x_j(m',n')}\\
&+\mu q^{-m_2n_2}\sum_{{-m_2=m_1+m+m'}\atop{-n_2=n_1+n+n'}}q^{n_1m'+nm_1+nm'}\p{\x_j(m,n)}\p{\x_i(m',n')}\\
=& 0
\end{align*}
as the first and the fourth canceled each other and the second one and
the third are negative to each other; and the last two kill each
other according to the calculation of the case
$[e_{1i}(m_1,n_1),e_{1i}(m_2,n_2)]$.

\begin{align*} &[e_{1i}(m_1,n_1),e_{ik}(m_2,n_2)]\\
=&q^{-m_1n_1}\mu q^{(-m_1-m_2)n_2}\p{\x_k(-m_1-m_2,-n_1-n_2)}\\
&-\sum
q^{n_1(m_2+m')+nm_1+n(m_2+m')+m'n_2}\\
&\cdot x_2(m_1+m_2+m+m',n_1+n_2+n+n')\p{\x_2(m,n)}\p{\x_k(m',n')}\\
&-\cdots\\
&-\sum
q^{n_1(m_2+m')+nm_1+n(m_2+m')+m'n_2}\\
&\cdot x_l(m_1+m_2+m+m',n_1+n_2+n+n')\p{\x_l(m,n)}\p{\x_k(m',n')}\\
&-\sum
q^{n_1m'+(n_2+n)m_1+(n_2+n)m'+mn_2}\\
&\cdot x_i(m_1+m_2+m+m',n_1+n_2+n+n')\p{\x_i(m',n')}\p{\x_k(m,n)}\\
&+\sum
q^{(m_1+m+m')n_2+n_1m'+nm_1+nm'}\\
&\cdot x_i(m_1+m_2+m+m',n_1+n_2+n+n')\p{\x_k(m,n)}\p{\x_i(m',n')}\\
=&q^{n_1m_2}e_{1k}(m_1+m_2,n_1+n_2).\end{align*}

For $i\neq k$:
\begin{align*} &[e_{1i}(m_1,n_1),e_{kk}(m_2,n_2)]\\
=&[-\sum
q^{n_1m'+nm_1+nm'}x_k(m_1+m+m',n_1+n+n')\p{\x_k(m,n)}\p{\x_i(m',n')},\\
&\sum q^{mn_2}x_k(m_2+m,n_2+n)\p{\x_k(m,n)}]\\
=& -\sum
q^{n_1m'+(n_2+n)m_1+(n_2+n)m'+mn_2}\\
&\quad\quad \cdot x_k(m_1+m_2+m+m',n_1+n_2+n+n')\p{\x_i(m',n')}\p{\x_k(m,n)}\\
&+\sum
q^{(m_1+m+m')n_2+n_1m'+nm_1+nm'}\\
&\quad\quad\cdot
x_k(m_1+m_2+m+m',n_1+n_2+n+n')\p{\x_k(m,n)}\p{\x_i(m',n')}\\
=&0.\end{align*}

 For $i\neq k$ and $k \neq j$:
\begin{align*} &[e_{1i}(m_1,n_1),e_{kj}(m_2,n_2)]\\
=&[-\sum
q^{n_1m'+nm_1+nm'}x_k(m_1+m+m',n_1+n+n')\p{\x_k(m,n)}\p{\x_i(m',n')}\\
&-\sum
q^{n_1m'+nm_1+nm'}x_j(m_1+m+m',n_1+n+n')\p{\x_j(m,n)}\p{\x_i(m',n')},\\
&\sum q^{mn_2}x_k(m_2+m,n_2+n)\p{\x_j(m,n)}]\\
=&-\sum
q^{n_1m'+(n_2+n)m_1+(n_2+n)m'+mn_2}\\
&\quad\quad\cdot x_k(m_1+m_2+m+m',n_1+n_2+n+n')\p{\x_i(m',n')}\p{\x_j(m,n)}\\
&+q^{(m_1+m+m')n_2+n_1m'+nm_1+nm'}\\
&\quad\quad\cdot x_k(m_1+m_2+m+m',n_1+n_2+n+n')\p{\x_j(m,n)}\p{\x_i(m',n')}\\
=&0.\end{align*}

\begin{align*} &[e_{1i}(m_1,n_1),e_{11}(m_2,n_2)]\\
=&[q^{-m_1n_1} \mu \p{\x_i(-m_1,-n_1},-\sum
q^{nm_2}x_k(m_2+m,n_2+n)\p{\x_k(m,n)}] \\
& +\sum_{k\neq i} [-\sum
q^{n_1m'+nm_1+nm'}x_k(m_1+m+m',n_1+n+n')\p{\x_k(m,n)}\p{\x_i(m',n')},\\
&-\sum q^{nm_2}x_k(m_2+m,n_2+n)\p{\x_k(m,n)}]\\
&+\sum_{k=2}^n [-\sum
q^{n_1m'+nm_1+nm'}x_k(m_1+m+m',n_1+n+n')\p{\x_k(m,n)}\p{\x_i(m',n')},\\
&-\sum q^{nm_2}x_i(m_2+m,n_2+n)\p{\x_i(m,n)}]\\
=&-q^{n_2m_1} \mu q^{-(m_1+m_2)(n_1+n_2)}
\p{\x_i(-(m_1+m_2),-(n_1+n_2))}
\\&+\sum_{k\neq i}\left(\sum
q^{n_1m'+(n_2+n)m_1+(n_2+n)m'+nm_2}\right.\\
&\quad\quad\cdot
x_k(m_1+m_2+m+m',n_1+n_2+n+n')\p{\x_k(m,n)}\p{\x_i(m',n')}\\
&-\sum
q^{(n_1+n+n')m_2+n_1m'+nm_1+nm')}\\
& \left.\quad\quad\cdot
x_k(m_1+m_2+m+m',n_1+n_2+n+n')\p{\x_k(m,n)}\p{\x_i(m',n')} \right)\\
&+\sum q^{n_1m'+(n_2+n)m_1+(n_2+n)m'+nm_2} \\
&\quad\quad\cdot
x_i(m_1+m_2+m+m',n_1+n_2+n+n')\p{\x_i(m,n)}\p{\x_i(m',n')}\\
&-\sum
q^{(n_1+n+n')m_2+n_1m'+nm_1+nm')}\\
& \quad\quad\cdot
x_i(m_1+m_2+m+m',n_1+n_2+n+n')\p{\x_i(m,n)}\p{\x_i(m',n')} \\
&+\sum_{k=2}^ n \sum
q^{n_1(m'+m_2)+nm_1+n(m'+m_2)+n'm_2}\\
&\quad\quad\cdot
x_k(m_1+m+m'+m_2,n_1+n+n'+n_2)\p{\x_k(m,n)}\p{\x_i(m',n')}\\
=& -q^{n_2m_1} \mu q^{-(m_1+m_2)(n_1+n_2)}
\p{\x_i(-(m_1+m_2),-(n_1+n_2))}\\
&+ q^{n_2m_1}\sum_{k\neq i}\sum
q^{n_1m'+nm_1+(n_2+n)m'+nm_2}\\
&\quad\quad\cdot
x_k(m_1+m_2+m+m',n_1+n_2+n+n')\p{\x_k(m,n)}\p{\x_i(m',n')}\\
\\ =&-q^{n_2m_1}e_{1i}(m_1+m_2,n_1+n_2).
\end{align*}

If $1<i,j,k,r\leqq n$,
\begin{align*} &[e_{ij}(m_1,n_1), e_{kr}(m_2,n_2)]\\
=&[\sum q^{mn_1}x_i(m_1+m,n_1+n)\p{\x_j(m,n)}, \sum
q^{mn_2}x_k(m_2+m,n_2+n)\p{\x_r(m,n)}]\\
=&\delta_{jk}\sum
q^{(m_2+m)n_1+mn_2}x_i(m_1+m_2+m,n_1+n_2+n)\p{\x_r(m,n)}\\
&-\delta_{ir}\sum q^{(m_1+m)n_2+mn_1}x_k(m_1+m_2+m,n_1+n_2+n)\p{\x_j(m,n)}\\
=&\delta_{jk} q^{m_2n_1}\sum
q^{m(n_1+n_2)}x_i(m_1+m_2+m,n_1+n_2+n)\p{\x_r(m,n)}\\
&-\delta_{ir}q^{m_1n_2}\sum
q^{m(n_1+n_2)}x_k(m_1+m_2+m,n_1+n_2+n)\p{\x_j(m,n)}\\
=&\delta_{jk} q^{m_2n_1} e_{ir}(m_1+m_2,n_1+n_2)-\delta_{ir}
q^{m_1n_2} e_{kj}(m_1+m_2,n_1+n_2).
\end{align*}

For $i\neq 1$,
\begin{align*} &[e_{ii}(m_1,n_1),e_{11}(m_2,n_2)]\\
=&[\sum q^{mn_1}x_i(m_1+m,n_1+n)\p{\x_i(m,n)},-\sum
q^{nm_2}x_i(m_2+m,n_2+n)\p{\x_i(m,n)}]\\
=&-\sum
q^{(m_2+m)n_1+nm_2}x_i(m_1+m_2+m,n_1+n_2+n)\p{\x_i(m,n)}\\
&+\sum
q^{(n_1+n)m_2+mn_1}x_i(m_1+m_2+m,n_1+n_2+n)\p{\x_i(m,n)}\\
=&0.\end{align*}

If $i\neq j$,
\begin{align*} &[e_{ij}(m_1,n_1),e_{11}(m_2,n_2)]\\
=&[\sum q^{mn_1}x_i(m_1+m,n_1+n)\p{\x_j(m,n)},\\
&-\sum q^{nm_2}x_i(m_2+m,n_2+n)\p{\x_i(m,n)}-\sum
q^{nm_2}x_j(m_2+m,n_2+n)\p{\x_j(m,n)}]\\
=&-\sum
q^{(m_2+m)n_1+nm_2}x_i(m_1+m_2+m,n_1+n_2+n)\p{\x_j(m,n)}\\
&+\sum
q^{(n_1+n)m_2+mn_1}x_i(m_1+m_2+m,n_1+n_2+n)\p{\x_j(m,n)}\\
=&0.\end{align*}

\begin{align*} &[e_{11}(m_1,n_1),e_{11}(m_2,n_2)]\\
=&\sum_{i=2}^l [-\sum q^{m_1n}x_i(m_1+m,n_1+n)\p{\x_i(m,n)},\\
&\quad \quad -\sum q^{m_2n}x_i(m_2+m,n_2+n)\p{\x_i(m,n)}]\\
=&\sum_{i=2}^l (\sum
q^{m_1(n_2+n)+m_2n}x_i(m_1+m_2+m,n_1+n_2+n)\p{\x_i(m,n)}\\
&-\sum
q^{m_2(n_1+n)+m_1n}x_i(m_1+m_2+m,n_1+n_2+n)\p{\x_i(m,n)})\\
=&-q^{m_1n_2}\mu\delta_{(m_1+m_2,n_1+n_2),(0,0)}\\
&+\sum_{i=2}^l\sum q^{m_1(n_2+n)+m_2n}x_i(m_1+m_2+m,n_1+n_2+n)\p{\x_i(m,n)}\\
&+q^{m_2n_1}\mu\delta_{(m_1+m_2,n_1+n_2),(0,0)}\\
&-\sum_{i=2}^l \sum q^{m_2(n_1+n)+m_1n}x_i(m_1+m_2+m,n_1+n_2+n)\p{\x_i(m,n)})\\
=&-q^{m_1n_2}e_{11}(m_1+m_2,n_1+n_2)+q^{m_2n_1}e_{11}(m_1+m_2,n_1+n_2).\end{align*}

For the proof of brackets involving $D_i$, we can refer to [Z] for
details. \qed\\

It is obvious that $V$ is a module of
$\widetilde{\frak{gl}_l(\bc_q)}$ generated by $1$, and
$\frak{n}_+.1=0$. Hence $V$ is the highest weight module of
$\widetilde{\frak{gl}_l(\bc_q)}$, and its highest weight vector is
$1$, with the weight $\alpha$, such that

$$\alpha(E_{11})=\mu, \alpha(E_{ii})=0, i\neq 1; \alpha(c_j)=\alpha(d_j)=0,
j=s,t$$

\begin{remark}The module given in [GZ1] is a special case of this module for $l=2$. The one in
[Z] is in a similar way of choosing different Borel subalgebra.
\end{remark}

\section{Irreducibility}
Since the centers $c_s, c_t$ act on $V$ trivially, we can view $V$ as a
module of $\frak{gl}_l\otimes \bc_q$. It is a weight module with the
Cartan subalgebra $\dot{\frak{h}}$, the highest weight is
$\alpha|_{\dot{\frak{h}}}$, we still denote it by $\alpha$.

We need to figure out the weight spaces of $V$. Before doing that, we
introduce some notations which will be used later.

Let
$$(M,N)=
\left(\begin{array}{ccc} (m_{21},n_{21}) & \cdots &
(m_{2k_2},n_{2k_2}) \\
(m_{31},n_{31}) & \cdots & (m_{3k_3},n_{3k_3}) \\
\cdots & \cdots & \cdots \\
(m_{l1},n_{l1}) & \cdots & (m_{lk_l},n_{lk_l})
\end{array} \right)
$$
with $k_i\in\bz_{\geq 0},(m_{ij},n_{ij})\in \bz^2$, and $m_{i1}\leq
m_{i2}\leq\cdots\leq m_{ik_i}$ for $i=2,\cdots l$. $(M,N)$ is a form
with $l-1$ lines, $k_i$ is the length of $(i-1)$-line, length is $0$
means this line is empty.

$(M,N)=(M',N')$ if two forms are the same; we call $(M,N)$ and
$(M',N')$ the same type if the length of every line in $(M,N)$ and
$(M',N')$ are the same. We denote the monomial $v_{(M,N)}$ as
$(M,N)$ defined above as

$$v_{(M,N)}=x_2(m_{21},n_{21})\cdots
x_2(m_{2k_2},n_{2k_2})\cdots x_l(m_{l1},n_{l1})\cdots
x_l(m_{lk_l},n_{lk_l})$$ If $k_i=0,i=2,\cdots,l$, this is $1$.

Since
$$e_{11}(0,0)(v_{(M,N)})=\mu-k_2-\cdots-k_l,
e_{ii}(0,0)(v)=k_i, i\neq 2,$$ the weight of $v_{(M,N)}$ is
$$\beta(E_{11})=\mu-k_2-\cdots-k_l, \beta(E_{ii})=k_i, i\neq 1.$$
Hence the weight space is spanned by all those monomials $v_{(M,N)}$
with $(M,N)$ the same type.

\begin{theorem} $V$ is irreducible module of $\frak{gl}_l\otimes \bc_q$
if and only if $\mu\neq 0$.
\end{theorem}

\noindent {\sl Proof}:\quad It is obvious that if $\mu=0$, $V$ is
reducible, with the largest submodule is the polynomials without
constant term.

 On the other hand, since this module is a highest weight module, we
only need to show that there is no other highest weight vector
except $1$ if $\mu\neq 0$.

Suppose $$v=\sum_{(M,N)} a_{(M,N)}v_{(M,N)}$$ be a highest weight
vector, here is a finite sum, and for all $(M,N)$ are the same
type, (i.e $v$ is in a weight space). Hence $\frak{n}_+$ acts as $0$ on
$v$.

At first we consider the case $j>i\geq 2$:

\begin{align*}
e_{ij}(m,n)\cdot v = & \sum_{(M,N)} a_{(M,N)}\sum_{r=1}^{k_j}
q^{m_{jr} n} x_i(m+m_{jr},n+n_{jr}) \cdot x_2(m_{21},n_{21}) \cdots
 \\
& \quad  x_2(m_{2k_2},n_{2k_2})\cdots \widehat{x_j(m_{jr},n_{jr})}
\cdots
x_l(m_{l1},n_{l1})\cdots x_l(m_{lk_l},n_{lk_l}) \\
=& \sum_{(M',N')}a'_{(M',N')}v_{(M',N')}=0,
\end{align*}
here $\widehat{ \quad\quad }$ means that this term is omitted (the
same for the following ). Suppose $(M,N)$ is not empty in some
$(j-1)$-th line, if we choose $(m,n)$ big enough, then in $(M',N')$,
the last one in the $(i-1)$-th line is $$(m_{jk_j}+m,n_{jk_j}+n),$$ so
if $a'_{(M',N')}=ca_{(M,N)}$ for some $(M,N)$, with $c$ is a nonzero
complex number, so $a'_{(M',N')}=0$, hence $a_{(M,N)}=0$.

Therefore for the highest weight vector, type of $(M,N)$ is empty in the
lines other than the first one, so we simplify
$(M,N)=((m_1,n_1),\cdots,(m_k,n_k))$ in the following.

Suppose that the highest weight vector is
$$v=\sum_{(M,N) \in K} a_{(M,N)} x_2(m_1,n_1) \cdots
x_2(m_k,n_k)$$ where $K$ is a set of some
$(M,N)=((m_1,n_1),\cdots,(m_k,n_k))$ with the same types.

In the following, we consider the action of $e_{12}(a,b)$. It is
easy to check that if $e_{12}(a,b)v_{(M,N)}$ and
$e_{12}(a,b)v_{(M',N')}$ are linearly dependent for all big enough
$a,b \in \bz$, then $v_{(M,N)}$ and $v_{(M',N')}$ must satisfy
$$
\sum_{i=1}^k m_i=\sum_{i=1}^k m'_i, \; \sum_{i=1}^k n_i=\sum_{i=1}^k
n'_i
$$

So we can assume that
$$
\sum_{i=1}^k m_i=m, \; \sum_{i=1}^k n_i=n, \; (M,N) =
((m_1,n_1),\cdots,(m_k,n_k)) \in K
$$
for some $m,n \in \bz$.

If $k=1$, then we have $v= a x_2(m,n)$, hence $$0=e_{12}(-m,-n)v =
aq^{-mn} \mu.$$ So we get $a=0$.

If $k=2$, then
$$
v=\sum_{i \in J} C_{m_i,n_i} x_2(m_i,n_i) x_2(m-m_i,n-n_i)
$$
where $J$ is an index set. Now let $I=\{(m_i,n_i) \,|\, i \in J\}$.
By symmetries, we can assume $m_i \leq \frac{m}{2}$ for $i \in J$
and if $(\frac{m}{2},n_j) \in I$ then $n_j \leq \frac{n}{2}$.

Now for $a,b \in \bz$, we have
$$
0=e_{12}(a,b) v = A x_2(m+a,n+b)
$$
where
\begin{align*}
A= & \mu q^{-ab} (C_{-a,-b}+C_{m+a,n+b}) \\
&+ \sum_{i \in J} C_{m_i,n_i} \left( q^{b(m-m_i)+an_i+n_i(m-m_i)} +
q^{bm_i +a(n-n_i) +m_i(n-n_i)} \right)
\end{align*}

First of all, for $b$ big enough and all $a \in \bz$, we have
\begin{align*}
0= & \sum_{i \in J} \left( C_{m_i,n_i}q^{an_i+n_i(m-m_i)}
q^{b(m-m_i)}+C_{m_i,n_i}q^{a(n-n_i)+m_i(n-n_i)}q^{bm_i}\right) \\
= & \sum_{i \in J'} A_{m_i} q^{b(m-m_i)} + \sum_{i \in J''} B_{m_i}
q^{bm_i}
\end{align*}
where $J' \subseteq J$ such that
$$\{m_i \,|\, i \in J' \} = \{m_i
\,|\, i \in J \} \text{ and } m_i \neq m_j \text{ if }i \neq j \in J';$$
 $J''
\subseteq J'$ such that
$$\{m_i \,|\, i \in J'' \} \dot\cup
\{\frac{m}{2}\} = \{m_i \,|\, i \in J' \}\text{ if }\frac{m}{2} \in \{m_i
\, |\, i \in J'\}$$ and $J'=J''$ otherwise; and where
$$
A_{m_i}=\sum_{(m_i,k) \in I} C_{m_i,n_i} \left( q^{a k+k(m-m_i)} +
\delta_{\frac{m}{2},m_i} q^{a(n-k) +m_i(n-k)} \right) \text{ for } i
\in J'
$$
and
$$
B_{m_i}=\sum_{(m_i,k) \in I} C_{m_i,n_i} q^{a (n-k)+(n-k)m_i} \text{
for } i \in J''.
$$

If $q$ is not a root of unity, then we have $q^{(m-m_i)}$, $i \in
J'$ and $q^{m_i}$, $i \in J''$ are $|J'|+|J''|$ distinct nonzero
complex numbers. Because of the nonvanishing of Vandermonde
dterminant, we have that $A_{m_i}=0$, $i \in J'$ and $B_{m_i}=0$, $i
\in J''$ hold for all $a \in \bz$. Now we have
$$
0=\sum_{i \in J} \left( C_{m_i,n_i}q^{an_i+n_i(m-m_i)} q^{b(m-m_i)}
+C_{m_i,n_i}q^{a(n-n_i)+m_i(n-n_i)}q^{bm_i}\right)
$$
for all $a,b \in \bz$.

If $q$ is a primitive $L-$th root of unity, then for $b$ big enough
and all $a \in \bz$, we have
\begin{align*}
0= & \sum_{i \in J}
\left( C_{m_i,n_i}q^{an_i+n_i(m-m_i)}
  q^{b(m-m_i)}+C_{m_i,n_i}q^{a(n-n_i)+m_i(n-n_i)}q^{bm_i}\right) \\
 = & \sum_{j=0}^{L-1} D_j q^{bj}
\end{align*}
where
$$
D_j=\sum_{i \in J'\atop j \equiv m-m_i (\text{mod }L)} A_{m_i} +
\sum_{i \in J''\atop j \equiv m_i (\text{mod }L)} B_{m_i}.
$$
Because of $q^i \neq q^j$ if $0 \leq i \neq j \leq L-1$, so we have
$D_j=0$ for all $0 \leq j \leq L-1$ and $a \in \bz$.  Hence
$$
0=\sum_{i \in J} \left( C_{m_i,n_i}q^{an_i+n_i(m-m_i)} q^{b(m-m_i)}
+C_{m_i,n_i}q^{a(n-n_i)+m_i(n-n_i)}q^{bm_i}\right)
$$
for all $a,b \in \bz$.

Now for any nonzero complex number $q$ and all $a,b \in \bz$, we
have
$$
A= \mu q^{-ab} (C_{-a,-b}+C_{m+a,n+b})=0.
$$

If $(\frac{m}{2},\frac{n}{2}) \in I$, taking $a=-\frac{m}{2}$ and
$b=-\frac{n}{2}$, we get $C_{\frac{m}{2},\frac{n}{2}}=0$.

If $(\frac{m}{2},\frac{n}{2}) \neq (m_i,n_i) \in I$, then
$(m-m_i,n-n_i) \notin I$.

Setting $a=-m_i$ and $b=-n_i$, we have
$$
0=\mu q^{-m_in_i} (C_{m_i,n_i}+C_{m-m_i,n-n_i})=\mu q^{-m_in_i}
C_{m_i,n_i},
$$
thus $C_{m_i,n_i}=0$. Now we know $C_{m_i,n_i}=0$ for all $i \in J$,
i.e., $v=0$.

If $k>2$ and $|K|=\alpha$, then
$$
v= \sum_{\beta =1}^\alpha C_\beta \, x_2(m_{\beta 1},n_{\beta 1})
x_2(m_{\beta 2},n_{\beta 2}) \cdots  x_2(m_{\beta k},n_{\beta k})
$$
where
$$C_\beta \in \bc, m_{\beta,k-1}+m_{\beta k} \geq m_{\gamma,
k-1} + m_{\gamma k} \text{ if } 1 \leq \beta < \gamma \leq \alpha$$
 and for
the same $m_{\beta,k-1}+m_{\beta k}$, giving the similar ordering on
$n_{\beta,k-1}+n_{\beta k}$. Now we rewrite $v$ as following
$$
v= \sum_{\beta =1}^\alpha C_\beta \, v_\beta x_2(m_{\beta,k-1}
,n_{\beta,k-1}) x_2(m_{\beta k},n_{\beta k}).
$$

Assume that
$$
v_1=v_2=\cdots =v_r \neq v_{r+1}
$$
and set
$$
m_{1,k-1}+m_{1,k}=\cdots=m_{r,k-1}+m_{r,k}=m'
$$
and
$$
n_{1,k-1}+n_{1,k}=\cdots=n_{r,k-1}+n_{r,k}=n'.
$$
Then
\begin{align*} v = & v_1 \sum_{\beta=1}^l C_\beta
x_2(m_{\beta,k-1} ,n_{\beta,k-1})
      x_2(m'-m_{\beta, k-1},n'-n_{\beta, k-1}) + \text{ other terms } \\
  = & v_1 \, v' + \text{ other terms }.
\end{align*}
Now we can easily see that $e_{12}(a,b)v'=0$ for big enough $b$ and
all $a \in \bz$. By above discussions, we have $v'=0$. Continuing this process,  we will have $v=0$.

Now, we must have $k=0$. i.e., $1$ is the only highest weight vector
up to a scalar multiple.   \qed

Together with the actions of $d_s$ and $d_t$ on $V$, we  can easily obtain
that

\begin{corollary}
$V$ is an irreducible module of $\widetilde{\frak{gl}_{l}(\bc_q)}$
if and only if $\mu\neq 0$.
\end{corollary}

If $\mu=0$, it is obvious that the module is reducible, with the
maximal submodule $W$ which consists of the polynomials without constant term. Here
we can get a parallel result of Thm 2.4 [GZ2].

\section*{Acknowledgment}
I am grateful to Professor Yun Gao for his encouragement and
 stimulating discussions. I would also like to thank Dr. Hongjia Chen for some helpful comments
on this work.

\end{document}